\documentstyle[amstex,amscd]{article}
\title{Counter-example to global Torelli problem for 
irreducible symplectic manifolds}
\author{Yoshinori Namikawa}
\date{ }
\begin{document}
\maketitle

A simply connected compact Kaehler manifold 
$X$ is an 
irreducible symplectic manifold if 
there is an everywhere non-degenerate 
holomorphic 2-form $\Omega$ on $X$ 
with $H^0(X, \Omega^2_X) = {\bold C}[\Omega]$. 
By definition, $X$ has even complex dimension. 
There is a canonical 
symmetric form $q_X$ on $H^2(X, {\bold Z})$, 
which is called the Beauville-Bogomolov form 
(cf. [Be]). 

On the other hand, since $X$ is Kaehler,  
$H^2(X, {\bold Z})$ has a natural Hodge structure 
of weight 2. If two irreducible symplectic 
manifolds $X$ and $Y$ are 
bimeromorphically equivalent, then 
there is a natural Hodge isometry between 
$(H^2(X, {\bold Z}), q_X)$ and 
$(H^2(Y, {\bold Z}), q_Y)$ (cf. [O, Proposition 
(1.6.2)], [Huy, Lemma 2.6]).  

Debarre [De] has constructed bimeromorphically 
equivalent irreducible symplectic manifolds $X$ and 
$Y$ such that $X$ and $Y$ are not isomorphic. 
This is a counter-example to the following 
problem. \vspace{0.2cm}

{\bf Biregular Torelli Problem}: 
{\em Let 
$X$ and $Y$ be irreducible symplectic manifolds 
of the same dimension 
such that there is a Hodge isometry $\phi: 
(H^2(X, {\bold Z}), q_X) \to (H^2(Y, {\bold Z}), q_Y)$.  
Is $X$ isomorphic to $Y$ ?} 
\vspace{0.12cm}

In this paper we shall give a counter-example 
to the following problems.
 \vspace{0.2cm}

{\bf Bimeromorphic Torelli Problem} 
(cf. [Mu, (5.10)], [Huy, (10.1)]): {\em Let 
$X$ and $Y$ be irreducible symplectic manifolds 
of the same dimension 
such that there is a Hodge isometry $\phi: 
(H^2(X, {\bold Z}), q_X) \to (H^2(Y, {\bold Z}), q_Y)$.  
Are $X$ and $Y$ bimeromorphically equivalent ?} 
\vspace{0.12cm} 

{\bf Polarized Torelli Problem}: 
{\em Let $(X, L)$ and $(Y, M)$ be polarized 
irreducible symplectic manifolds of the same 
dimension such that there is a Hodge 
isometry $\phi: (H^2(X, {\bold Z}), q_X) 
\to (H^2(Y, {\bold Z}), q_Y)$ with 
$\phi([L]) = [M]$. Is $(X, L)$ isomorphic 
to $(Y, M)$ as a polarized variety ?} 
\vspace{0.12cm}   

In the paragraphs {\bf 1}, ..., {\bf 4} we 
 construct a counter-example to 
Bimeromorphic Torelli Problem, and in {\bf 5} 
we discuss Polarized Torelli Problem. 
\vspace{0.12cm}

{\bf 1}. Let $T$ be a complex torus of dimension 2 and 
let ${\mathrm Hilb}^{n+1}(T)$ be the Hilbert 
scheme (or Douady space) that parametrizes 
length $n+1$ points on $T$. There is a Hilbert-
Chow map $h: 
{\mathrm Hilb}^{n+1}(T) \to {\mathrm Sym}^{n+1}(T)$. 
Here ${\mathrm Sym}^{n+1}(T)$ is the $n+1$ symmetric 
product of $T$. Given a group structure on $T$, we 
have a holomorphic map $\alpha : 
{\mathrm Sym}^{n+1}(T) \to T$ 
by sending $[p_1, ..., p_{n+1}]$ to $\Sigma p_i$. 
Let $K^n(T)$ be the fiber over $0 \in T$ 
of the composite of two maps:    
$$ {\mathrm Hilb}^{n+1}(T) \to T. $$ 
Then $K^n(T)$ becomes an irreducible symplectic 
manifold of dim $2n$ (cf. [Be]). Let $q$ be the 
Beauville-Bogomolov form on $K^n(T)$.  

Let us consider the case where $n = 2$. 
Put $\bar{K}^2(T) := \alpha^{-1}(0)$. Then 
$\bar{K}^2(T)$ has only quotient singularities. 
The singular locus $\Sigma$ 
of $\bar{K}^2(T)$ is isomorphic to $T$. 
There is  a bimeromorphic map 
$h_0 : K^2(T) \to \bar{K}^2(T)$ 
and its exceptional locus $E$ is an irreducible 
divisor of $K^2(T)$. 
The general fiber of the map 
$E \to \Sigma$ is isomorphic to ${\bold P}^1$. 
Let $F$ be a resolution of $E$. Then there 
is a holomorphic surjective map from 
$F$ to $\Sigma$. This map coincides with the 
Albanese map of $F$. Therefore, the Albanese 
variety of $F$ is isomorphic to $T$.  
\vspace{0.2cm} 

{\bf 2}. By [Yo, Lemma(4.10), Proposition(4.11)], there 
is a natural Hodge isometry for $n \geq 2$: 

$$ H^2(K^{n}(T), {\bold Z}) \cong 
H^2(T, {\bold Z})\oplus {\bold Z}\delta, $$ 
where the left hand side is equipped with the 
Beauville-Bogomolov form $q$ and the right hand 
side is the direct sum of two lattices   
$(H^2(T, {\bold Z}), ( , ))$ and ${\bold Z}\delta$ 
with $\delta^2 = -2(n+1)$. The Hodge structure 
on the right hand side is given by 

$ H^{2,0} := H^{2,0}(T), $ 

$ H^{1,1} := H^{1,1}(T) \oplus {\bold C}\delta,$ and 

$ H^{0,2} := H^{0,2}(T)$.   

By the construction (cf. [Yo, (4.3.1)]) we 
have $2\delta 
= [E]$. 
\vspace{0.2cm}

{\bf 3}. Let $T$ be a complex torus of dimension 2 such that 

(1) the dual torus $T^*$ of $T$ is not isomorphic 
to $T$, and   

(2) the Neron-Severi group of $T$ is 
trivial: $\mathrm{NS}(T) = 0$.  \vspace{0.2cm}

Since $H_1(T, {\bold Z}) \cong H^1(T^*, {\bold Z})$, 
we have a pairing $H^1(T, {\bold Z}) \times 
H^1(T^*, {\bold Z}) \to {\bold Z}$. This induces 
the non-degenerate pairing $\wedge^2 H^1(T, {\bold Z}) 
\times \wedge^2 H^1(T^*, {\bold Z}) \to {\bold Z}$. 
Since $\wedge^2 H^1(T, {\bold Z}) = 
H^2(T, {\bold Z})$ and $\wedge^2 H^1(T^*, {\bold Z}) 
= H^2(T^*, {\bold Z})$, we have the non-degenerate 
pairing $H^2(T, {\bold Z}) \times H^2(T^*, {\bold Z}) 
\to {\bold Z}$. This pairing induces an 
isomorphism $H^2(T^*, {\bold Z}) \cong 
\mathrm{Hom}(H^2(T, {\bold Z}), {\bold Z})$. 
By the cup-product $H^2(T, {\bold Z}) \times 
H^2(T, {\bold Z}) \to H^4(T, {\bold Z}) \cong 
{\bold Z}$,\footnote{The isomorphism 
$H^4(T, {\bold Z}) \cong {\bold Z}$ is given by 
the orientation defined by the complex structure}  
$\mathrm{Hom}(H^2(T, {\bold Z}), {\bold Z})$ 
is identified with $H^2(T, {\bold Z})$. 
Therefore we obtain a canonical isomorphism 
$$ \alpha_T : H^2(T^*, {\bold Z}) 
\cong H^2(T, {\bold Z}). $$
By Shioda [Sh], this isomorphism is a Hodge 
isometry. 
By the paragraph 
{\bf 2}, there is a Hodge isometry  
$H^2(K^2(T^*), {\bold Z}) \to H^2(K^2(T), {\bold Z})$ 
extending this Hodge isometry.  
Now we shall prove the following. \vspace{0.2cm}

{\bf Proposition}. {\em There are no 
bimeromorphic maps from 
$K^2(T)$ to $K^2(T^*)$.} \vspace{0.15cm}

{\bf 4}. ({\em Proof of Proposition}): We put 
$X := K^2(T)$ and $Y := K^2(T^*)$. 
Assume that there is a bimeromorphic map 
$f : X --\to Y$. Since $\omega_X \cong 
{\cal O}_X$ and $\omega_Y \cong {\cal O}_Y$ 
we see that $f$ is an isomorphism in 
codimension 1. Therefore, $f$ induces an 
isomorphism $f^*: H^2(Y, {\bold Z}) \cong 
H^2(X, {\bold Z})$. Moreover, this map 
induces an isomorphism 
${\mathrm Pic}(Y) \cong {\mathrm Pic}(X)$. 
By {\bf 2} and the assumption (2) of {\bf 3} 
we see that ${\mathrm Pic}(X) = {\bold Z}\delta$ 
and ${\mathrm Pic}(Y) = {\bold Z}\delta^*$. 
Since $f^*$ is a Hodge isometry with respect 
to Beauville-Bogomolov forms (cf. [O, Proposition 
(1.6.2)], [Huy, Lemma 2.6]), 
we conclude that $f^*(\delta^*) 
= \delta$ or $f^*(\delta^*) = -\delta$. 
But, since $2\delta$ (resp. $2\delta^*$) 
is represented by $E$ (resp. $E^*$)(see 
paragraph {\bf 2}), the latter 
case does not occur because $X$ and $Y$ 
are Kaehler manifolds. Now, since 
$f$ is an isomorphism in codimension 1, 
$f$ induces a bimeromorphic map between 
$E$ and $E^*$. As in {\bf 1}, let $F$ 
(resp. $F^*$) be a resolution of $E$ 
(resp. $E^*$). Since $F$ and $F^*$ are 
bimeromorphic, there should be a natural 
isomorphism between their Albanese 
varieties. By {\bf 1}, ${\mathrm Alb}(F) = T$ 
and ${\mathrm Alb}(F^*) = T^*$.   
This contradicts the assumption (1) of 
{\bf 3}. \vspace{0.2cm}

{\bf Remark 1}. As is well known, three Torelli 
problems are affirmative for $K3$ surfaces. 
Our counter-example is valid only for $K^n(T)$ 
with $n \geq 2$. The situation is 
quite different for the Kummer 
surface $K^1(T)$. 
First note that the exceptional locus of 
the bimeromorphic map $h_0 : K^1(T) 
\to \bar{K}^1(T)$ consists of sixteen 
$(-2)$-curves $C_i$. Let $M$ be the 
smallest primitive sublattice of 
$H^2(K^1(T), {\bold Z})$ containing 
$\oplus {\bold Z}[C_i]$. 
   
In $H^2(K^1(T), {\bold Z})$, the 
primitive sublattice 
$H^2(\bar{K}^1(T), {\bold Z})$ is 
isomorphic to $(H^2(T, {\bold Z}), 2< , >)$, 
where $< , >$ is the cup product on $T$. 
Since $H^2(K^1(T), {\bold Z})$ is a 
unimodular lattice, $H^2(\bar{K}^1(T), 
{\bold Z}) \oplus M$ is of finite index 
$> 1$ in $H^2(K^1(T), {\bold Z})$.     

As in {\bf 3}, let $T^*$ be the 
dual torus of $T$. Then the Hodge 
isometry $H^2(T^*, {\bold Z}) 
\cong H^2(T, {\bold Z})$ 
constructed in {\bf 3} 
induces an isometry between 
$H^2(\bar{K}^1(T^*), {\bold Z}) 
\oplus M^*$ and 
$H^2(\bar{K}^1(T), {\bold Z}) 
\oplus M$. 
However, this isometry 
does not extend to an isometry 
between $H^2(K^1(T^*), {\bold Z})$ 
and $H^2(K^1(T), {\bold Z})$.  
\vspace{0.2cm} 
 
{\bf 5}. One can 
construct an example in the 
category of projective varieties. 
For example, let  
$T$ be an Abelian surface  
such that 

(1) the dual torus $T^*$ of $T$ is not isomorphic 
to $T$, and   

(2) $\mathrm{NS}(T) = {\bold Z}[H]$ with 
$H^2 = 6$. 

Put $H^* := {\alpha_T}^{-1}(H)$. 
One can check that $H^*$ is an 
ample class of $H^2(T^*, {\bold Z})$. 
Then, in the proof of 
Proposition, ${\mathrm Pic}(X) 
= {\bold Z}H \oplus {\bold Z}\delta$ 
and ${\mathrm Pic}(Y) = {\bold Z}H^* 
\oplus {\bold Z}\delta^*$. Since 
$f^*$ induces an isometry between  
these two lattices, it is easily 
checked that $f^*(H^*) = H$ or $-H$, 
and $f^*(\delta^*) = \delta$ or 
$-\delta$. By the same argument as 
Proposition, the latter cases are 
excluded and we have $f^*(H^*) = H$ 
and $f^*(\delta^*) = \delta$. 
Therefore, we 
conclude that $K^2(T)$ and $K^2(T^*)$ 
are not bimeromorphically equivalent.    

Let us identify $H^2(K^2(T), {\bold Z})$ 
(resp. $H^2(K^2(T^*), {\bold Z})$) 
with $H^2(T, {\bold Z}) \oplus {\bold Z}\delta$ 
(resp. $H^2(T^*, {\bold Z}) \oplus {\bold Z}
\delta^*$) as in {\bf 2}.    
If $m > 0$ is a sufficiently large 
integer, then $[L] := m[H] - \delta \in 
H^2(K^2(T), {\bold Z})$ and 
$[M] := m[H^*] - \delta^* \in 
H^2(K^2(T^*), {\bold Z})$ are 
both ample classes. The Hodge isometry 
 $H^2(K^2(T), {\bold Z}) 
\to H^2(K^2(T^*), {\bold Z})$ in {\bf 3} sends 
$[L]$ to $[M]$. So this gives a counter-example 
to Polarized Torelli Problem.
\vspace{0.12cm}

{\bf Remark 2}. If we replace the Abelian 
surface $T$ in {\bf 5} by the one with 
$H^2 = 4$, then a Fourier-Mukai transform 
induces an isomorphism 
$\phi : K^2(T^*) \to K^2(T)$ such that 
$\phi^*(H) = 5H^* - 4\delta^*$, 
$\phi^*(\delta) = 6H^* - 5\delta^*$ 
(cf. [Yo, Propositions (3.5),(4.9)]). 
\vspace{0.12cm}

{\bf Remark 3}. Let $T$ be an Abelian surface 
and $T^*$ its dual. Put $X = K^2(T)$ and 
$Y = K^2(T^*)$. 
Let $D(X)$ (resp. $D(Y)$) be the bounded 
derived category of coherent sheaves on 
$X$ (resp. $Y$). Then there is an equivalence 
of categories between $D(X)$ and $D(Y)$. 
In fact, the symmetric group $G := S_3$ acts 
on $T^3 := T \times T \times T$ by the 
permutation. $G$ acts also on $N := 
\{(x,y,z) \in T^3; x + y + z = 0 \}$.  
Let $G$-${\mathrm Hilb}(N)$ be the 
$G$-Hilbert scheme of the $G$-variety 
$N$ (cf. [B-K-R]). Then the 
irreducible component of 
$G$-${\mathrm Hilb}(N)$ containing the free 
orbits becomes $X$ (cf. [Ha, Theorem 
5.1]). Now apply [B-K-R, Corollary 1.3] 
to the diagram 
$$  X \to N/G \gets N. $$ 
Then we have an equivalence of categories 
$\psi: D(X) \cong D^G(N)$, where 
$D^G(N)$ is the bounded derived category of coherent 
G-sheaves on $N$. 
Similarly we get an equivalence of 
categories $\psi^* : D(Y) \cong 
D^G(N^*)$, where $N^* := \{(x,y,z) \in 
(T^*)^3; x + y + z = 0 \}$. Since 
$N^*$ is the dual Abelian variety of 
$N$, the Fourier-Mukai transform induces 
an equivalence between $D^G(N)$ and 
$D^G(N^*)$. Therefore, we have an equivalence 
between $D(X)$ and $D(Y)$.  
\vspace{0.12cm}

{\bf Question}: {\em Let $X$ and $Y$ be 
two irreducible symplectic manifolds 
such that there is a Hodge isometry 
$(H^2(X, {\bold Z}), q_X) \cong 
(H^2(Y, {\bold Z}), q_Y)$. Then, 
is there an equivalence of categories 
between $D(X)$ and $D(Y)$?}  
\vspace{0.12cm}

When $X$ and $Y$ are $K3$ surfaces, 
$D(X)$ and $D(Y)$ are equivalent if 
and only if there is a Hodge isometry 
between transcendental lattices of 
$H^2(X, {\bold Z})$ and $H^2(Y, {\bold Z})$ 
([Or]).  
\vspace{0.12cm} 

There are another series of examples of 
irreducible symplectic manifolds, namely 
those which are deformation equivalent to 
${\mathrm Hilb}^n(S)$ with $S$ being a $K3$ 
surface. I do not know any negative evidence 
for Bimeromorphic Torelli Problem for 
such manifolds.        
   
\vspace{0.2cm}

{\em Acknowledgement}. The author 
thanks A. Fujiki, S. Mukai, Y. Kawamata
and D. Huybrechts for discussion or  
comments on this paper.

\vspace{0.2cm}

{\bf Erratum.}(added in Aug. 4, 2002): 
 Remark 3 is not correct.  
In fact, the $G$ action on $N$ induces 
a natural $G$ action on $N^*$, but this 
action {\em differs} from the permutaion 
of $(x,y,z) \in N^*$. Let $D^G(N^*)$ be the 
same as in Remark 3 and let $D^G(N^*)'$ 
be the bounded derived category of coherent 
$G$-sheaves with respect to this induced $G$-action. 
There is an equivalence between $D^G(N)$ 
and $D^G(N^*)'$ by the same reason as Remark 3, 
but there is no such equivalence between 
$D^G(N)$ and $D^G(N^*)$.     
The quotient variety $N^*/G$ for the induced 
$G$-action on $N^*$ is actually a symplectic 
V-manifold, but it has no crepant resolutions.   
\vspace{0.2cm}

\begin{center}
Department of Mathematics, 
Graduate school of science, 
Osaka University, Toyonaka, Osaka 560, Japan 
\end{center}  
 
\end{document}